\documentclass[12pt]{article}
\usepackage{amssymb,latexsym}
\makeatletter
\newfont{\footsc}{cmcsc10 at 8truept}
\newfont{\footbf}{cmbx10 at 8truept}
\newfont{\footrm}{cmr10 at 10truept}
\makeatother
\title{A Local Characterization of Combinatorial Multihedrality in Tilings
 }
 
\author{Nikolai Dolbilin \thanks{Supported, in part, by RFBR grant 05-01-00170
 and SSS 2185.2003.1.}
\\
\small Steklov Mathematical Institute\\[-0.8ex]
\small Gubkin 8\\[-0.8ex]
\small Moscow 117966, Russia\\[-0.8ex]
\small \texttt{dolbilin@mi.ras.ru}
\\[3mm]
\small and \\[3mm]
Egon Schulte\thanks{Supported, in part, by NSA-grants
H98230-04-1-0116
and H98230-05-1-0027.} \\
\small Northeastern University\\[-0.8ex]
\small Department of Mathematics\\[-0.8ex]
\small Boston, MA 02115, USA\\[-0.8ex]
\small \texttt{schulte@neu.edu}}

\date{Version as of \today}

\newenvironment{proof}{\par \noindent {\em
Proof}:\quad}{\hspace*{\fill} $\Box$ \par \vspace*{1ex}}
\newenvironment{proofthm}{\par \noindent {\em
Proof of Theorem~\ref{chper}}:\quad}{\hspace*{\fill} $\Box$ \par \vspace*{1ex}}
\makeatletter
\@addtoreset{equation}{section}
\makeatother

\newtheorem{theorem}{Theorem}[section]
\newtheorem{lemma}[theorem]{Lemma}
\newtheorem{conjecture}[theorem]{Conjecture}
\newtheorem{problem}[theorem]{Problem}
\newtheorem{corollary}[theorem]{Corollary}
\newtheorem{definition}[theorem]{Definition}
\newtheorem{remark}[theorem]{Remark}
\newcommand{\bpf}{\begin{proof}}
\newcommand{\epf}{\end{proof}}
\newcommand{\blem}{\begin{lemma}}
\newcommand{\bthm}{\begin{theorem}}
\newcommand{\bcon}{\begin{conjecture}}
\newcommand{\bpro}{\begin{problem}}
\newcommand{\bcor}{\begin{corollary}}
\newcommand{\bdefin}{\begin{definition}}
\newcommand{\brem}{\begin{remark}}
\newcommand{\elem}{\end{lemma}}
\newcommand{\ethm}{\end{theorem}}
\newcommand{\econ}{\end{conjecture}}
\newcommand{\epro}{\end{problem}}
\newcommand{\ecor}{\end{corollary}}
\newcommand{\edefin}{\end{definition}}
\newcommand{\erem}{\end{remark}}
\newcommand{\beq}{\begin{equation}}
\newcommand{\bqy}{\begin{eqnarray*}}
\newcommand{\bry}{\begin{array}}
\newcommand{\eeq}{\end{equation}}
\newcommand{\eqy}{\end{eqnarray*}}
\newcommand{\ery}{\end{array}}
\newcommand{\eqref}[1]{(\ref{#1})}

\newcommand{\Ga}{{\mit \Gamma}}

\newcommand{\CT}{{\cal T}}

\newcommand{\CC}{{\cal C}}

\newcommand{\GaT}{\Gamma({\cal T})}
\newcommand{\BE}{\mathbb{E}}

\newcommand{\BH}{\mathbb{H}}
\begin{document}
\maketitle

\begin{abstract}
\noindent
A locally finite face-to-face tiling of euclidean $d$-space by convex polytopes is called {\em combinatorially multihedral\/} if its combinatorial automorphism group has only finitely many orbits on the tiles. The paper describes a local characterization of combinatorially multihedral tilings in terms of centered coronas. This generalizes the Local Theorem for Monotypic Tilings, established in \cite{dolsch}, which characterizes the case of combinatorial tile-transitivity.
\end{abstract}

\section{Introduction}
\label{intro}
Given a geometric structure in space it is generally a difficult task to formulate local conditions that predetermine its global properties. In the present paper, these structures are locally finite face-to-face tilings of Euclidean $d$-space $\BE^d$ by convex $d$-polytopes. In our earlier paper~\cite{dolsch}, we obtained a local characterization of combinatorial tile-transitivity of monotypic tilings (with a single combinatorial prototile) in $\BE^d$ by convex $d$-polytopes;  the result is the ``Local Theorem for Monotypic Tilings". This characterization is expressed in terms of combinatorial symmetry properties of large enough neighborhood complexes (centered coronas) of tiles. The theorem is a combinatorial analogue of the ``Local Theorem for Tilings", which describes a local characterization of isohedrality (geometric tile-transitivity) of monohedral tilings (with a single isometric prototile) in $\BE^d$ (see \cite{ddsg,ds}). Both theorems are closely related to the ``Local Theorem for Delone Sets" of \cite{ddsg}, which locally characterizes those sets among the uniformly discrete sets in $\BE^d$ which are orbits under a crystallographic group (see also 
\cite{dls}).

This paper characterizes combinatorial multihedrality of locally finite face-to-face tilings of $\BE^d$ by convex $d$-polytopes. The local conditions employed are strong enough to force the combinatorial automorphism group of the tiling to have only finitely many orbits on the tiles. In some sense, this characterization can be viewed as attempting to map out the fuzzy border between combinatorial ``periodicity" and ``non-periodicity" of tilings using transitivity properties of groups that act on them (see also \cite{es1,es2}). The new ``Local Theorem" is a combinatorial analogue of a local characterization of periodic (or crystallographic) tilings in $\BE^d$ established in \cite{dol2}; the latter characterization, in turn, is a tiling version of a similar such characterization for crystallographic point sets obtained in \cite{dolshto}. The combinatorial analogue we establish here, is considerably harder to prove than its geometric counterpart (see Section~\ref{cry}) and depends in an essential way on the simply-connectedness of the underlying space, unlike the original version.

\section{Some terminology}
\label{termin}

Recall that a {\em tiling\/} $\CT$ of euclidean $d$-space $\BE^d$ is a countable family of closed topological $d$-balls of $\BE^d$, the {\em tiles} of $\CT$,  which cover $\BE^{d}$ without gaps and overlaps (see, for instance, Gr\"unbaum \& Shephard~\cite{gs}). In this paper, the tiles will always be convex $d$-polytopes (see \cite{gcp}). (For the present paper it would actually suffice to require the tiles to be homeomorphic images of convex polytopes. However, as in our previous paper \cite{dolsch}, to avoid unnecessary technicalities which blur the essence of the overall argument, it is appropriate here to assume convexity of the tiles.) All tilings are taken to be {\em locally finite\/}, meaning that each point of $\BE^d$ has a neighborhood that meets only finitely many tiles.  Moreover, all tilings (by convex polytopes) are assumed to be {\em face-to-face\/}, that is, the intersection of any two tiles is a face of each tile, possibly the empty face. For a face-to-face tiling $\CT$, the set of all faces of the tiles (with $\CT$ itself adjoined as improper face of rank $d+1$), ordered by inclusion, becomes a lattice, called {\em the face-lattice\/} of $\CT$.

A locally finite face-to-face tiling $\CT$ of $\BE^d$  by convex polytopes is said to be
{\em combinatorially multihedral\/} (or {\em combinatorially crystallographic\/}) if its combinatorial automorphism group $\GaT$ has only finitely many orbits on the tiles. (Since the tiles are convex polytopes, there is a group of homeomorphisms of $\BE^d$ which is isomorphic to $\GaT$ and has the same action on the face-lattice of $\CT$ as $\GaT$. Hence we could have defined combinatorial multihedrality in terms of this group of topological automorphisms. However, here we will always work directly with $\GaT$.)  If there are exactly $n$ tile orbits under $\GaT$, we call $\CT$ {\em combinatorially tile-$n$-transitive\/}. (Note that combinatorially tile-$n$-transitive tilings which are normal in the sense of \cite{gs}, are $n$-homeohedral in the sense of \cite{gs}, and vice versa.)  When $n=1$, the tiling $\CT$ is {\em combinatorially tile-transitive\/}. 

These notions are straightforward combinatorial analogues of similar such notions involving the geometric symmetry group $G(\CT)$ of $\CT$. Recall that a tiling $\CT$ of $\BE^d$ is {\em periodic\/} (or {\em geometrically crystallographic\/}) if and only if $\CT$ has only finitely many orbits of tiles under $G(\CT)$, or, equivalently, $G(\CT)$ is a crystallographic group (a discrete group of isometries of $\BE^d$ with compact fundamental domain). If there are precisely $n$ orbits under $G(\CT)$, then $\CT$ is said to be $n$-{\em isohedral\/}, or simply {\em isohedral\/} if $n=1$ (see \cite{gs}).

Any two tiles $P$ and $Q$ of a locally finite face-to-face tiling $\CT$ of $\BE^d$ can be joined by a finite sequence of tiles 
\begin{equation}
\label{seqdist}
P=P_0,P_1,\ldots,P_{n-1},P_{n}=Q 
\end{equation} 
of $\CT$ such that $P_{j-1}\cap P_{j}$ is a face of $P_{j-1}$ and $P_{j}$ of dimension at least $d-2$, for $j=1,\ldots,n$; we call $n$ the {\em length\/} of the sequence. The minimum length of a sequence joining tiles $P$ and $Q$ as in \eqref{seqdist} is called the {\em distance\/} of $P$ and $Q$ in $\CT$ and is denoted by $d(P,Q)$. In proofs we often employ sequences \eqref{seqdist}, in which the tiles $P_{j-1}$ and $P_j$ share a facet ($(d-1)$-face) for $j=1,\ldots,n$; any two tiles $P$ and $Q$ of $\CT$ can be joined by such a sequence. 

There are variants of this distance function which require any two consecutive tiles in a sequence to intersect in a face of dimension at least $l$, for some fixed $l$ with $0 \leq l \leq d-1$. The present distance function corresponds to the case $l=d-2$. We shall explain in a moment why we distinguished the number $d-2$. 

The local characterization of combinatorial multihedrality established in this paper is based on the notion of a centered corona of $\CT$. Let $P$ be a tile of $\CT$, and let $k \geq 0$ be an integer.  The {\em $k^{th}$ corona of $P$\/}, denoted by $\CC_k(P)$, is the subcomplex of $\CT$ consisting of the tiles $Q$ of $\CT$ with $d(P,Q) \leq k$, and their faces. In particular, the $0^{th}$ corona $\CC_0(P)$ is the face-lattice of $P$, and, if $k\geq 1$, the $k^{th}$ corona $\CC_k(P)$ is the set of faces of tiles that intersect a tile in $\CC_{k-1}(P)$ in a face of dimension at least $d-2$. By definition, a corona is a complex, not a set of tiles or a union of tiles (this differs from the use of the term in \cite{ds,enone}). Two distinct tiles $P$ and $Q$ can have the same $k^{th}$ corona for some $k$, that is, $\CC_k(P) = \CC_k(Q)$ and hence $\CC_j(P) = \CC_j(Q)$ for each $j\geq k$ (see \cite[Figure 1]{dolsch} for an example with $k=1$). Therefore we shall distinguish coronas by their tile of reference.  A {\em centered $k^{th}$ corona\/} is a pair $(P,C_k(P))$ consisting of a tile $P$ of $\CT$, the {\em center\/} of the centered $k^{th}$ corona, and its $k^{th}$ corona $\CC_{k}(P)$ in $\CT$.  When it is clear from the context we simply denote $(P,\CC_k(P))$ by $\CC_k(P)$. 

Note that the notion of centered corona depends on the distance function employed. For the type of distance functions which require in their definition that any two consecutive tiles in a sequence intersect in a face of dimension at least $l$, for some fixed $l$ with $0 \leq l \leq d-1$, the size of the centered coronas is smallest (and hence most desirable) when $l$ is largest. However, the statement of our Theorem~\ref{chper} fails for centered coronas based on the distance function with $l=d-1$, but holds for those with $l\leq d-2$. Thus $l=d-2$ is the natural choice, yielding the smallest centered coronas for which the theorem is valid.

Two centered $k^{th}$ coronas $\CC_k(P)$ and $\CC_k(P')$ of a tiling $\CT$ are said to be {\em isomorphic\/} if there exists an isomorphism of complexes $\alpha: \CC_k(P) \longrightarrow \CC_k(P')$ with $\alpha (P)=P'$; such a map $\alpha$ is called an {\em isomorphism of centered $k^{th}$ coronas\/}. We will only consider isomorphisms between coronas which map center to center. In particular, for a centered $k^{th}$ coronas $\CC_k(P)$, we denote by $\Ga(\CC_k(P))$ its {\em automorphism group\/}, that is, the stabilizer of the center $P$ in the full automorphism group of the (non-centered) corona $\CC_k(P)$. 

An automorphism of the centered $(k-1)^{st}$ corona $\CC_{k-1}(P)$ of a tile $P$ may not in general be extendable to an automorphism of the centered $k^{th}$ corona $\CC_k(P)$. However, if it is, then it does extend uniquely. (Note here that every isomorphism between a centered corona and a subcomplex of $\CT$ is uniquely determined by its effect on a single flag of the corona~\cite[Lemma 2.1]{dolsch}.)  On the other hand, the restriction of each automorphism of $\CC_k(P)$ to $\CC_{k-1}(P)$ is also an automorphism of $\CC_{k-1}(P)$, which, in turn, uniquely determines the original automorphism of $\CC_k(P)$.  Therefore,
\beq
\label{subgroup}
\Gamma (\CC_{k-1}(P))\supseteq \Gamma (\CC_{k}(P)) \quad (k\geq 1),
\eeq
with equality occurring if and only if each automorphism of $\CC_{k-1}(P)$ is extendable to an automorphism of $\CC_k(P)$. Moreover, each automorphism of $\CT$ that fixes $P$ restricts to an automorphism of $\CC_k(P)$ for each $k\geq 0$, and is uniquely determined by this restriction. Thus, for a tile $P$ of $\CT$, we have an infinite chain of subgroup relationships,
\beq
\label{chain}
\Ga(P) = \Ga(\CC_0(P)) \supseteq \Ga(\CC_1(P)) \supseteq \ldots \supseteq \Ga(\CC_k(P)) \supseteq \ldots \supseteq
\Ga_{P}(\CT ), 
\eeq
with the combinatorial automorphism group $\Ga(P)$ of $P$ on the left and the stabilizer $\Ga_{P}(\CT)$ of $P$ in $\Ga(\CT)$ on the right. Since $\Ga(P)$ is a finite group, the number of proper descents in \eqref{chain} is finite and is bounded by the number of prime divisors of $|\Ga(P)|$ (counted with multiplicity).

Note that, if the centered $k^{th}$ coronas $\CC_{k}(P)$ and $\CC_{k}(P')$ of two tiles $P$ and $P'$ are isomorphic under an isomorphism of centered coronas $\alpha$, then $\Ga(\CC_k(P))$ and $\Ga(\CC_k(P'))$ are isomorphic under the group isomorphism 
$\beta \rightarrow \alpha\beta\alpha^{-1}$.

Once again, let $\CT$ be a locally finite face-to-face tiling of $\BE^d$.  For $k\geq 0$, let 
$N_{k} := N_{k}(\CT)$ denote the number of isomorphism classes of centered $k^{th}$ coronas of tiles in $\CT$; this number may be infinite. If $N_k$ is finite for some $k\geq 0$, then the number $N_0$ of combinatorial isomorphism classes of tiles of $\CT$ must be finite as well; in fact, otherwise the $k^{th}$ centered coronas of the tiles in an infinite family of mutually non-isomorphic tiles would yield an infinite family of mutually non-isomorphic centered $k^{th}$ coronas. However, the converse may not be true; that is, a finite $N_0$ may not imply a finite $N_k$. (Note here that the tiling may not be normal in the sense of Gr\"unbaum \& Shephard~\cite{gs}.) From now on we always assume that $N_k$ is finite for each $k\geq 0$.

If $P$ and $P'$ are tiles of $\CT$ with isomorphic centered $k^{th}$ coronas $\CC_{k}(P)$ and $\CC_{k}(P')$, then their centered $(k-1)^{st}$ coronas $\CC_{k-1}(P)$ and $\CC_{k-1}(P')$ are also isomorphic. Hence
\beq
\label{nk}
N_{k-1} \leq N_{k} \quad (k\geq 1) . 
\eeq
In other words, there are at least as many isomorphism classes of centered $k^{th}$ coronas as there are isomorphism classes of centered $(k-1)^{st}$ coronas. The following lemma characterizes the occurrence of the equality sign in \eqref{nk}.

\blem
\label{remone}
Suppose $k$ is a positive integer with $N_{k-1} = N_{k} < \infty$. Then the centered 
$k^{th}$ coronas of two tiles of $\CT$ are isomorphic if and only if their centered $(k-1)^{st}$ coronas are isomorphic.
\elem

\bpf
Let $P$ and $P'$ be tiles of $\CT$. We show that isomorphism of $\CC_{k-1}(P)$ and $\CC_{k-1}(P')$ implies isomorphism of $\CC_{k}(P)$ and $\CC_{k}(P')$. In fact, by our assumption, if $P_{1},\ldots,P_{N_{k}}$ are tiles whose centered $k^{th}$ coronas form a full set of representatives for the isomorphism classes of centered $k^{th}$ coronas, then the corresponding centered $(k-1)^{st}$ coronas must form a full set of representatives for the isomorphism classes of centered $(k-1)^{st}$ coronas (no two of the latter can be isomorphic). Now we can argue as follows. We know that
$\CC_{k}(P) \simeq \CC_{k}(P_{i})$ and $\CC_{k}(P') \simeq \CC_{k}(P_{j})$ for some $i$ and $j$, and hence $\CC_{k-1}(P) \simeq \CC_{k-1}(P_{i})$ and $\CC_{k-1}(P') \simeq \CC_{k-1}(P_{j})$. Therefore, if $\CC_{k-1}(P) \simeq \CC_{k-1}(P')$, then necessarily $i=j$ and hence $\CC_{k}(P) \simeq \CC_{k}(P')$.   
\epf

As before, let $P_{1},\ldots,P_{N_{k}}$ be tiles whose centered $k^{th}$ coronas form a full set of representatives for the isomorphism classes of centered $k^{th}$ coronas of $\CT$. For $i=1,\ldots,N_{k}$ and $j \leq k$, let 
\[ \Gamma_{j}^{i} := \Gamma(\CC_{j}(P_{i})) ; \]
that is, $\Gamma_{j}^{i}$ is the automorphism group of the centered $j^{th}$ corona of the $i^{th}$ tile $P_{i}$. Then, by \eqref{subgroup},
\beq
\label{gk}
\Gamma_{k-1}^{i} \supseteq \Gamma_{k}^{i} \quad (i = 1,\ldots,N_{k}) . 
\eeq
Note that, for each $i$, the $i^{th}$ isomorphism class of centered $k^{th}$ coronas determines the corresponding groups $\Gamma_{j}^{i}$ ($j\leq k)$ up to group isomorphism (induced by conjugation of maps). In particular, if $P$ is a tile whose centered $k^{th}$ corona represents the $i^{th}$ isomorphism class of centered $k^{th}$ coronas of $\CT$, then any two such groups $\Gamma_{j}^{i}$ and $\Gamma_{j'}^{i}$ coincide if and only if the corresponding groups 
$\Gamma(\CC_{j}(P))$ and $\Gamma(\CC_{j'}(P))$ coincide.

\section{A generalization of the Local Theorem for Monotypic Tilings}
\label{newloc}

The following Local Theorem for combinatorially multihedral tilings is a combinatorial analogue of the Local Theorem for crystallographic tilings obtained in \cite{dol2}, and a generalization of the Local Theorem for Monotypic Tilings obtained in \cite{dolsch}.

\bthm
\label{chper}
Let $\CT$ be a locally finite face-to-face tiling of $\BE^d$ by convex $d$-polytopes. Then $\CT$ is combinatorially multihedral if and only if for some positive integer $k$ the following properties hold:
\begin{enumerate}
\item $N_{k}(\CT)$ is finite and $N_{k-1}(\CT) = N_{k}(\CT)$.
\item $\Gamma_{k-1}^{i} = \Gamma_k^i\; (i = 1,\ldots ,N_{k}(\CT))$, where $\Gamma_{k-1}^{i}$ and $\Gamma_k^{i}$ represent the automorphism groups of the centered coronas at levels $k-1$ and $k$, respectively, of a tile associated with the $i^{th}$ isomorphism class of centered $k^{th}$ coronas of $\CT$.
\end{enumerate}
In particular, in this case, if $P$ is a tile of $\CT$, then $\Gamma(\CC_{k}(P)) = \Gamma_{P}(\CT)$. Moreover, if $n$ is a positive integer, then $\CT$ is combinatorially tile-$n$-transitive if and only if for some positive integer $k$ the two properties hold with $N_{k}(\CT)=n$.
\ethm

The proof follows a similar line of arguments as the proof for the Local Theorem of Monotypic Tilings in \cite{dolsch}. The main challenge is the sufficiency part. This is based on a key lemma 
(Lemma~\ref{lem6} below), which in turn follows from Lemma~\ref{lem1} and a series of lemmas established in \cite{dolsch}.  

The statement of Lemma~\ref{lem1} is basically identical with \cite[Lemma 3.2]{dolsch}, but the proof given in \cite{dolsch} is no longer valid and needs to be adapted to the present situation. Throughout, bear in mind that $k$ is a positive integer satisfying the two conditions of Theorem~\ref{chper}. 

\blem
\label{lem1}
Let $P$ and $P'$ be tiles of $\CT$, and let $\bar{\alpha}: \CC_{k-1}(P)\rightarrow \CC_{k-1}(P')$ be an isomorphism of centered $(k-1)^{st}$ coronas. Then $\bar{\alpha}$ extends uniquely to an isomorphism of centered $k^{th}$ coronas $\alpha: \CC_{k}(P)\rightarrow \CC_{k}(P')$. 
\elem

\bpf
Most arguments of the old proof go through again step by step. First, every automorphism of the centered $(k-1)^{st}$ corona $\CC_{k-1}(P)$ extends uniquely to an automorphism of the centered $k^{th}$ corona $\CC_{k}(P)$; this basically follows as before, now appealing to the second condition $\Gamma_{k-1}^{i} = \Gamma_k^i$ of the theorem, where $i$ is such that $\CC_{k}(P)$ represents the $i^{th}$ isomorphism class of centered $k^{th}$ coronas. The next step in the proof then requires the existence of any isomorphism whatsoever between the centered $k^{th}$ coronas $\CC_{k}(P)$ and $\CC_{k}(P')$. In \cite{dolsch} this follows directly by assumption, but not so here. However, in the present context the existence of an isomorphism is still implied by the first condition $N_{k-1}=N_k$ of the theorem, bearing in mind Lemma~\ref{remone}; in fact, by assumption, $\CC_{k-1}(P)$ and $\CC_{k-1}(P')$ are isomorphic (under $\bar{\alpha}$), and hence $\CC_{k}(P)$ and $\CC_{k}(P')$ are isomorphic as well. The remaining steps of the old proof then carry over to the present situation.
\epf
\smallskip

The next Lemma~\ref{lem6} says that every local isomorphism $\alpha$ between centered $k^{th}$ coronas extends uniquely to an automorphism of the whole tiling $\CT$, thereby becoming a global isomorphism. Its proof is quite involved and depends heavily on the simply-connectedness of the underlying space. The key step in the proof is the extension of $\alpha$ along sequences of tiles as in \eqref{seqdist} and their centered coronas. 

\blem
\label{lem6}
Let $P$ and $P'$ be tiles of $\CT$, let $\alpha: \CC_{k}(P) \rightarrow \CC_{k}(P')$ be an isomorphism of centered $k^{th}$ coronas. Then there exists a combinatorial automorphism $\varphi$ of $\CT$ (that is, $\varphi\in \GaT$) which extends $\alpha$, that is, $\varphi |_{\CC_{k}(P)} = \alpha$. Moreover, $\varphi$ is uniquely determined.
\elem

\bpf
In essence this is \cite[Lemma 3.7]{dolsch}, whose proof is based on \cite[Lemmas 3.2-3.6]{dolsch}. However, in the present context the old proof requires modification as follows. First, \cite[Lemma 3.2]{dolsch} is replaced by the new Lemma~\ref{lem1} above. Then, with this replacement, the statements and lengthy proofs of \cite[Lemmas 3.3-3.6]{dolsch} remain valid and establish, together with Lemma~\ref{lem1}, the existence and uniqueness of $\varphi$.
\epf

We now have the tools to prove Theorem~\ref{chper}.
\medskip

\begin{proofthm}
We first obtain the necessity of the two conditions. Let $\CT$ be combinatorially multihedral, and in particular let $\CT$ be combinatorially tile-$n$-transitive. If $P$ and $P'$ are tiles of $\CT$ equivalent under $\GaT$, then every automorphism of $\CT$ mapping $P$ to $P'$ induces an isomorphism between the centered $k^{th}$ coronas of $P$ and $P'$, for each integer $k\geq 0$. Hence, by inequality \eqref{nk}, $(N_{k})_{k\geq 0}$ is a non-decreasing sequence of positive integers bounded by $n$. Then the first condition of the theorem clearly holds for some positive integer $k$ (in fact, for all large enough $k$).  Moreover, we can choose $k$ such that the second condition of the theorem is satisfied as well. In fact, since the automorphism group $\Ga (P)$ of each tile $P$ is finite, there can only be a finite number of proper descents in the subgroup chain \eqref{chain} associated with $P$, so each pair $\Gamma(\CC_{k-1}(P))$, $\Gamma(\CC_{k}(P))$ of consecutive subgroups in this chain must coincide when $k$ is sufficiently large. Therefore, since there are only finitely many combinatorially different tiles in $\CT$, the second condition holds for large enough $k$. This establishes necessity of the two conditions. We later see that we also have $N_{k} = n$ in this case. 

Next we prove sufficiency of the two conditions. Suppose that $k$ is a positive integer satisfying both conditions of the theorem. We shall prove that any two tiles of $\CT$ associated with the same isomorphism class of centered $k^{th}$ coronas are in fact equivalent under $\GaT$. In fact, let $P$ and $P'$ be two tiles of $\CT$ whose centered $k^{th}$ coronas are isomorphic, and let $\alpha :\CC_k(P)\longrightarrow \CC_k(P')$ be an isomorphism of centered $k^{th}$ coronas. Then, by Lemma~\ref{lem6}, $\alpha$ admits an extension to a combinatorial automorphism $\varphi$ of $\CT$ which maps $P$ to $P'$, so in particular $P$ and $P'$ are equivalent under $\GaT$. It follows that any two tiles associated with the same isomorphism class of centered $k^{th}$ coronas are equivalent under $\GaT$.  Hence $\GaT$ has at most $N_{k}$ orbits on the tiles, so in particular $\CT$ is combinatorially multihedral. Notice that the number of orbits of $\GaT$ on the tiles must be precisely $N_{k}$. In fact, any two tiles equivalent under $\GaT$ must have isomorphic centered $k^{th}$ coronas, so fewer orbits under $\GaT$ would lead to fewer isomorphism class of centered $k^{th}$ coronas in $\CT$.  Moreover, the above argument applied with $P=P'$ shows that $\Gamma(\CC_{k}(P)) = \Gamma_{P}(\CT)$. In fact, in the above we have $\varphi\in\Gamma_{P}(\CT)$; this proves one inclusion, bearing in mind that $\varphi$ extends the element $\alpha$ of $\Gamma(\CC_{k}(P))$, while the other follows from \eqref{chain}. 

Finally, we just established that, under the two conditions of the theorem on the parameter $k$, the number of orbits of $\GaT$ on the tiles equals the number $N_{k}$ of centered $k^{th}$ coronas of $\CT$.  In particular, this also settles the claim at the end of the necessity proof saying that $N_{k} = n$.
\end{proofthm}

We remark that the second condition of Theorem~\ref{chper} cannot be ignored in general. Examples (already with $n=1$) in the plane have been discussed in \cite[p.6]{dolsch} and are based on results in Gr\"unbaum-Shephard~\cite{gsinc}. 

The following corollary deals with a special case when the second condition of Theorem~\ref{chper} is trivially satisfied.  Recall that a convex $d$-polytope $P$ in $\BE^d$ is {\em combinatorially asymmetric\/} if its combinatorial automorphism group $\Ga(P)$ is trivial. 

\bcor
\label{corol}
Let $\CT$ be a locally finite face-to-face tiling of $\BE^d$ by combinatorially asymmetric convex $d$-polytopes, and let $n$ be a positive integer. Then $\CT$ is combinatorially multihedral with $n$ tile orbits under $\GaT$ if and only if $N_{k-1}(\CT) = N_{k}(\CT)=n$ for some positive integer~$k$.
\ecor

\brem
\label{remnew}
As mentioned earlier there are variants of the distance function which require in their definition that any two consecutive tiles in \eqref{seqdist} intersect in a face of dimension at least $l$, for some fixed $l$ with $0\leq l \leq d-1$. Each variant leads to a new notion of centered $k^{th}$ corona, and then, in turn, to a variant of Theorem~\ref{chper} based on this notion (and proved in a similar way), for each $l$ with $l \leq d-2$. The case $l=d-2$ yields the original theorem itself. However, if $l=d-1$, the analogous result already fails in the tile-transitive case (see \cite[Remark 3.9]{dolsch}); in fact, the proof of Lemma~\ref{lem6} requires $l\leq d-2$. 
\erem

As pointed before, Theorem~\ref{chper} also generalizes to locally finite face-to-face tilings of $\BE^d$ by {\em topological\/} $d$-polytopes (homeomorphic images of convex $d$-polytopes), and hence to such tilings of hyperbolic $d$-space $\BH^d$ as well; that is, convexity of the tiles is not really needed.

\section{Crystallographic tilings}
\label{cry}

Theorem~\ref{chper} is a combinatorial analogue of the following Theorem~\ref{geomchper} from \cite{dol2} (see also \cite{dolshto}), which characterizes periodic (or geometrically crystallographic) tilings in $\BE^d$. On the surface, the two theorems look quite
similar. However, there are very significant differences; in particular, neither one implies the
other. Theorem~\ref{geomchper} involves global isometries of the ambient space and pairwise
congruence classes of centered tile coronas, whereas Theorem~\ref{chper} is expressed in terms of (local) combinatorial isomorphisms of subcomplexes (namely, centered coronas) and combinatorial isomorphism classes of such subcomplexes. Moreover, more importantly, the proof of Theorem~\ref{chper} relies (in Lemma~\ref{lem6}) heavily on the simply-connectedness of the underlying space, while the proof of Theorem~\ref{geomchper} does not. Notice also that the two concepts of coronas are different.

We briefly review some definitions. Once again, let $\CT$ be a locally finite face-to-face tiling of $\BE^d$. The {\em $k^{th}$ tile corona\/} of a tile $P$ in $\CT$, denoted by $C_k(P)$, is the set of tiles $Q$ of $\CT$ with $d(P,Q) \leq k$, where the distance function $d(.\,,.)$ is as before. Then the old $k^{th}$ corona $\CC_{k}(P)$ is simply the complex of all faces of tiles in the $k^{th}$ tile corona $C_{k}(P)$. As before, the {\em centered $k^{th}$ tile corona\/} of $P$ is the pair $(P,C_k(P))$ consisting of $P$, the {\em center\/} of the centered $k^{th}$ tile corona (usually dropped from the notation), and the $k^{th}$ tile corona $C_{k}(P)$. Two centered $k^{th}$ tile coronas $C_k(P)$ and $C_k(P')$ of $\CT$ are {\em pairwise congruent\/} if there exists an isometry $\alpha$ of $\BE^d$ with $\alpha(P) = P'$ and $\alpha(C_{k}(P)) = C_{k}(P')$; that is, $\alpha$ induces a bijection between the sets $C_{k}(P)$ and $C_{k}(P')$ which maps $P$ to $P'$. Thus we can speak about the {\em pairwise congruence class\/} of a $k^{th}$ tile corona.

Let $P$ be a tile of $\CT$, and let $G(P)$ be its geometric symmetry group in $\BE^d$. The 
{\em symmetry group\/} $G_{k}(P)$ of the centered $k^{th}$ tile corona $C_{k}(P)$ is defined by
\[ G_{k}(P) := \{ \alpha \in G(P) \mid \alpha(C_{k}(P)) = C_{k}(P) \}. \]
Its elements are the {\em symmetries\/} of $C_{k}(P)$ (in general they are not symmetries of $\CT$). Then $G_{k}(P)$ is a subgroup of $G(P)$ for each $k$, and $G_{0}(P)=G(P)$. Moreover, we have the following infinite chain of subgroup relationships,
\beq
\label{geomchain}
G(P) = G_{0}(P) \supseteq G_{1}(P) \supseteq \ldots \supseteq
G_{k}(P) \supseteq \ldots \supseteq
G_{P}(\CT ), 
\eeq
with $G(P)$ on the left and the stabilizer $G_{P}(\CT)$ of $P$ in the geometric symmetry group $G(\CT)$ on the right. In fact, if $k \geq 1$, then every symmetry of $C_k(P)$ preserves (combinatorial) distance from $P$ and hence is a symmetry of $C_{k-1}(P)$; moreover, being an isometry, it is uniquely determined by its effect on $C_{k-1}(P)$. Similarly, if $k\geq 0$, then every symmetry of $\CT$ that fixes $P$ is a symmetry of $C_k(P)$. Since $G(P)$ is a finite group, there are at most a finite number of proper descents in \eqref{geomchain}; this number is bounded by the number of prime divisors of $|G(P)|$ (counted with multiplicity). 

For $k\geq 0$, let $M_{k}:=M_{k}(\CT)$ denote the number of pairwise congruence classes of centered $k^{th}$ tile coronas of tiles in $\CT$. Then $M_0$ is the number of congruence classes of tiles in $\CT$. If $P$ and $P'$ are tiles of $\CT$ with pairwise congruent centered $k^{th}$ tile coronas $C_{k}(P)$ and $C_{k}(P')$, then their centered $(k-1)^{st}$ coronas $C_{k-1}(P)$ and $C_{k-1}(P')$ are also pairwise congruent. Hence 
\beq
\label{mk}
M_{k-1} \leq M_{k} \quad (k\geq 1) , 
\eeq
and therefore $M_0$ is finite if $M_k$ is finite (for at least one $k$). However, in the present context (of isometric congruence), the converse is true as well. In fact, $M_0$ is finite if and only if $M_k$ is finite for each $k$. Now suppose that $M_0$ is finite. Let $k\geq 0$, and let $P_{1},\ldots,P_{M_{k}}$ be tiles whose centered $k^{th}$ tile coronas form a full set of representatives for the pairwise congruence classes of centered $k^{th}$ tile coronas of $\CT$. For $i=1,\ldots,M_{k}$ and 
$j \leq k$, let $G_{j}^{i} := G_{j}(P_{i})$; that is, $G_{j}^{i}$ is the symmetry group of the centered $j^{th}$ tile corona of $P_{i}$. Then, by \eqref{geomchain},
\beq
\label{geomgk}
G_{k-1}^{i} \supseteq G_{k}^{i}  \quad (i = 1,\ldots,M_{k}) . 
\eeq
Note that the $i^{th}$ pairwise congruence class of centered $k^{th}$ tile coronas determines the corresponding groups $G_{j}^{i}$ up to conjugacy. 

Now periodicity of tilings in $\BE^d$ can be characterized by the following Local Theorem.

\bthm
\label{geomchper}
Let $\CT$ be a locally finite face-to-face tiling of $\BE^d$ by convex $d$-polytopes, and let $\CT$ have only finitely many congruence classes of tiles. Then $\CT$ is periodic if and only if there exists a positive integer $k$ with the following properties:
\begin{enumerate}
\item $M_{k-1}(\CT) = M_{k}(\CT)$.
\item $G_{k-1}^{i} = G_{k}^{i} \; (i = 1,\ldots,M_{k}(\CT))$, where $G_{k-1}^{i}$ and $G_{k}^{i}$ 
represent the symmetry groups of the centered tile coronas at levels $k-1$ and $k$, respectively, of a tile associated with the $i^{th}$ pairwise congruence class of centered $k^{th}$ tile coronas of $\CT$. 
\end{enumerate}
In particular, in this case, if $P$ is a tile of $\CT$, then $G_{k}(P) = G_{P}(\CT)$. Moreover, if $n$ is a positive integer, then $\CT$ is $n$-isohedral if and only if for some positive integer $k$ the two properties hold with $M_{k}(\CT)=n$.
\ethm

The comments about generalizations made in Remark~\ref{remnew} also apply in the present context. In fact, Theorem~\ref{geomchper} even holds in the previously excluded case $l=d-1$. Yet another variant of the present theorem applies to arbitrary (not necessarily face-to-face) tilings; here the distance function is such that two tiles are at distance $1$ if and only if they intersect.

\end{document}